\documentclass[12pt]{article}
\usepackage[T2A]{fontenc}
\usepackage[utf8]{inputenc}
\usepackage{amssymb,amsfonts,amsmath,color}
\usepackage{amssymb,amsmath,amsfonts,amsthm,amscd,float}
\usepackage{latexsym,indentfirst,verbatim,mathtext,cite,enumerate,geometry}
\def\Span{\mathrm{Span}}
\def\Aut{\mathrm{Aut}}
\def\Der{\mathrm{Der}}
\def\charr{\mathrm{char}}
\def\RB{\textrm{RB}}
\def\Imm{\textrm{Im}}
\def\rb{\textrm{rb}}
\def\id{\textrm{id}}
\def\ext{\textrm{ext}}

\begin{document}

\sloppy


\begin{center}
{\Large Burde's series of simple pre-Lie algebras}

\smallskip

Vsevolod Gubarev
\end{center}

\begin{abstract}
We describe automorphisms, derivations, and Rota---Baxter operators
on the series of simple pre-Lie algebras found by D. Burde in 1998.

\end{abstract}

\section{Introduction}

Pre-Lie algebras appeared during 1960s independently in affine geometry
(E.B. Vinberg~\cite{Vinberg}; J.-L. Koszul~\cite{Koszul})
and ring theory (M. Gerstenhaber~\cite{Gerstenhaber}).
Pre-Lie algebras are also known under names of Vinberg algebras, Koszul algebras,
left- or right-symmetric algebras (LSAs or RSAs), Gerstenhaber ones.
An algebra~$A$ is a pre-Lie algebra, if
\begin{equation}\label{preLie}
(x\cdot y)\cdot z - x\cdot (y\cdot z)
 = (y\cdot x)\cdot z - y\cdot (x\cdot z)
\end{equation}
holds for all $x,y,z\in A$.

Note that every pre-Lie algebra $A$ under the commutator 
$[a,b] = ab-ba$ is a Lie algebra. 
Thus, pre-Lie algebras are Lie-admissible algebras.

The problem of classification of simple finite-dimensional 
pre-Lie algebras is open and it seems to be very hard. 
For example, in~\cite{Pozhidaev}, it was shown 
that starting with a~finite-dimensional pre-Lie algebra~$A$ 
with a~nonzero product, one can construct 
a~simple finite-dimensional pre-Lie algebra, 
which contains~$A$ as a~subalgebra.
In 1996, D.~Burde showed that all simple pre-Lie algebras 
over $\mathrm{gl}_n(\mathbb{C})$ may be obtained as deformations 
of the matrix algebra.

In~1998, D.~Burde found~\cite{Burde} a series of simple pre-Lie algebras~$I_n$.
A.~Mizuhara based on the studies of H. Shima~\cite{Shima} gave 
a~list of simple pre-Lie algebras over a solvable Lie algebra~\cite{Mizuhara}, 
including~$I_n$.
In~\cite{Burde}, it was proved that $I_2$ is the only 
simple pre-Lie algebra of dimension~2 over~$\mathbb{C}$.
However, already in dimension~3, there are an infinite number 
of pairwise non-isomorphic simple pre-Lie algebras over~$\mathbb{C}$~\cite{Burde}.
All simple pre-Lie algebras obtained from~$I_4$ by infinitesimal deformations 
were described by D.~Burde in~\cite{Burde}.
Simple real forms on $I_n$ were given in~\cite{Kong}.

In~\cite{MajidTao}, it was shown that the quantum special linear algebra 
$SL_q(2)$ over $\mathbb{C}$ with real~$q$ 
may be quantised from the pre-Lie product on~$\mathrm{su}_2^*$
which occurs to be isomorphic to~$I_3$ (over~$\mathbb{R}$).

E.B. Vinberg showed~\cite{Vinberg} the one-to-one correspondence
between homogeneous regular convex domains in $V = \mathbb{R}^n$ 
and so called clans, pre-Lie algebras $(V,\cdot)$ endowed with 
a~positive linear form such that the operator $L_x$ 
of left multiplication on~$x$ has for all $x\in V$ 
only real eigenvalues.
Such clans admit a~principal idempotent~$e$; 
if the eigenspace of the operator $R_e$ 
of right multiplication on~$e$ 
corresponding to~1 is one-dimensional, 
then the clan is called elementary. 
In~\cite[Proposition\,10.7]{Shima2}, H. Shima proved that an~elementary clan 
is isomorphic to~$I_n$. 
Some close result was recently obtained by D. Fox~\cite[Lemma\,7.2]{Fox}.

In~\cite{Pozhidaev2}, simple finite-dimensional pre-Lie algebras from
the variety of algebras satisfying the identities 
$[[x,y],[z,t]] = [x,y]([z,t]u) = 0$ were studied.
In~\cite{Chapoton}, some simple infinite-dimensional 
graded pre-Lie algebras of special form were classified.

Pre-Lie algebras which additionally satisfy the identity $(xy)z = (xz)y$ 
are called Novikov ones~\cite{Balinskii}. 
It is known~\cite{Zelmanov} that over a field of characteristic zero, 
every simple finite-dimensional Novikov algebra is isomorphic to a~field. 
In~\cite{Filippov,Osborn,Xu}, simple infinite-dimensional 
Novikov algebras were constructed, and in~\cite{Osborn} 
derivations of some of them were described.

In~\cite{AnBai,Bai,preLie}, classification of derivations and Rota---Baxter operators on $I_2$ as well as description of Rota---Baxter operators of weight~0 on $I_3$ were obtained.
In the current work, we describe all main operators on $I_n$: automorphisms, derivations, and Rota---Baxter operators of an arbitrary weight~$\lambda$.

More detailed, we prove that $\Aut(I_n)\cong O_{n-1}(F)$ and $\Der(I_n)\cong \mathfrak{o}_{n-1}(F)$ (Theorem~1).
Further, we show that every Rota---Baxter operator~$R$ of weight~$\lambda$ on $I_n$ satisfies the equation $R^2+\lambda R = 0$ (Theorem~2).
In particular, it means that every Rota---Baxter operator of nonzero weight on $I_n$ corresponds to a direct decomposition of $I_n$ into a sum of two subalgebras.
We show that one of the subagebras has to be Lagrangian and another one is a~span of a complement and~$e_n$.
Finally, we discuss for $I_n$ simplicity of its anticommutator algebra~$I_n^{(+)}$ as well as its infinite-dimensional analogue~$I_\infty$ (Remarks~2 and~3).

Throughout the work, we suppose that characteristic 
of the ground field~$F$ does not equal~2.

\section{Preliminaries}

Let us give important examples of pre-Lie algebras and then introduce the simple pre-Lie algebra~$I_n$.

{\bf Example 1}~\cite{Svinolupov}.
Let $V$ be a finite-dimensional vector space over $\mathbb{R}$ 
with the scalar product $(\cdot,\cdot)$ 
and let $(0\neq)a\in V$. Then $V$ under the product
$u\circ v = (u,v)a + (u,a)v$
is a~simple pre-Lie algebra.

{\bf Example 2}~\cite{Burde2}.
Denote by $\partial_i$ the derivation $\frac{\partial}{\partial x_i}$ 
in the algebra $A = F[x_1,\ldots,x_n]$, where $F$ is a~field.
Let us define on the space 
$$
\mathrm{Der}(A) 
 = \left\{ \sum\limits_{i=1}^n f_i \partial_i \mid f_i\in A \right\}
$$
a new product 
$u\partial_i \circ v\partial_j = u\partial_i(v)\partial_j$.
The obtained algebra $(\mathrm{Der}(A),\circ)$ is a pre-Lie algebra,
which under the commutator
$[a,b] = a\circ b - b\circ a$ is exactly Witt algebra. 

Let $U_n(F)$ denote the associative algebra 
of non-stricly upper-triangular matrices. 
A.Z. Anan'in studied in~\cite{Anan'in} 
the varities of algebras which may be embedded into $U_n(F)$
under the ordinary product.
Our main object~$I_n$ is also connected to the algebra $U_n(F)$.

{\bf Example 3}~\cite{Burde,Nijenhuis}.
Define a linear map $\tau\colon M_n(F)\to U_n(F)$ as follows,
$$
\tau(e_{ij}) = \begin{cases}
e_{ij}, & i<j, \\
0, & i>j, \\
e_{ii}/2, & i = j.
\end{cases}
$$ 
Then $U_n(F)$ is a pre-Lie algebra under the product 
$x\circ y = xy + \tau(xy^T + yx^T)$, where $T$ denotes the transposition.

Although $(U_n(F),\circ)$ is not simple,
its ideal $I_n = \Span\{e_{1j}\mid j=1,\ldots,n\}$ is. 
Denote $e_k = e_{1\,n+1-k}$, $k=1,\ldots,n$, the linear basis of $I_n$, 
and define the product~$\cdot$ by the formulas
$$
e_n \cdot e_n = 2e_n,\quad
e_n \cdot e_j = e_j,\quad
e_j \cdot e_j = e_n,\quad j=1,\ldots,n-1,
$$
all omitted products are assumed equal to zero.
Hence, $I_n$ is a~simple pre-Lie algebra~\cite{Burde}.
Note that $I_n$ is simple over a field of any characteristic, 
not necessarily of characteristic zero, as it was considered in~\cite{Burde}. 

When $n = 1$, $I_1$ is isomorphic to the field~$F$.
For $n\geq2$, $I_n$ is not only non-associative, moreover it is 
not power-associative:
$$
e_1 = (e_1\cdot e_1)\cdot e_1
 \neq e_1\cdot (e_1\cdot e_1) = 0.
$$
The pre-Lie algebra~$I_n$ is not complete when $\mathrm{char}\,F\neq0$, 
i.\,e. $\mathrm{tr}(R_x)$ is not zero for all $x\in I_n$. We have $\mathrm{tr}(R_{e_n}) = 2$.

Denote $\Span\{e_1,\ldots,e_{n-1}\}$ by $J_n$.
Every element $x\in I_n$ may be presented as 
$x = v + \alpha e_n$, where $v\in J_n$ and $\alpha\in F$.
Define a symmetric bilinear form on the space~$J_n$ as follows: 
for $a = a_1 e_1 + \ldots + a_{n-1}e_{n-1}$ and 
$b = b_1 e_1 + \ldots + b_{n-1}e_{n-1}$, we put
$(a,b) = a_1 b_1 + \ldots + a_{n-1}b_{n-1}$.

For 
$x = v + \alpha e_n$ and 
$y = u + \beta e_{n}$, where $v,u\in J_n$, $\alpha,\beta\in F$, 
we have by the multiplication table
\begin{equation} \label{In:MultTable}
xy
 = \alpha u
 + ((v,u)+2\alpha\beta)e_n.
\end{equation}

In~\cite{Kong}, it was noted that if we deal over~$\mathbb{R}$ and 
we have a scalar product on $I_n$ such that $e_1,\ldots,e_n$ 
is an orthonormal basis, then the product on~$I_n$ is defined by Example~1 with $a = e_n$.
Indeed,
$$
x\circ y 
 = (x,e_n)y + (x,y)e_n
 = \alpha(u+\beta e_n) + ((v,u)+\alpha\beta)e_n
 = xy.
$$

\section{Automorphisms and derivations}

For an arbitrary field~$F$, the orthogonal group $O_n(F)$ and the 
orthogonal Lie algebra $\mathfrak{o}_n(F)$ are defined as follows,
$$
O_n(F) = \{A\in M_n(F)\mid AA^T = A^TA = E\}, \quad
\mathfrak{o}_n(F) = \{A\in M_n(F)\mid A^T = -A\},
$$
here $E$ denotes the identity matrix.

In~\cite{Bai}, it was shown there are no nonzero derivations on~$I_2$.
Let us describe the automorphism group and the Lie algebra of derivations on $I_n$.

{\bf Theorem 1}. 
Let $n>1$ and $\charr F\neq2$. We have 
$\Aut(I_n)\cong O_{n-1}(F)$ and $\Der(I_n)\cong \mathfrak{o}_{n-1}(F)$.

{\sc Proof}.
Let $\varphi\in \Aut(I_n)$, write
$\varphi(e_i) = v_i + \alpha_i e_n$, where $v_i\in J_n$ and $\alpha_i\in F$.
Thus, 
$$
\varphi(e_i e_j)
 = \varphi(e_i)\varphi(e_j)
 = \alpha_i v_j + ( (v_i,v_j) + 2\alpha_i \alpha_j )e_n.
$$
This equality is equivalent to the system of relations 
\begin{gather}
\alpha_i v_i = v_n,\quad (v_i,v_i) + 2\alpha_i^2 
 = \alpha_n,\quad i<n; \label{Aut:ii} \\
\alpha_n v_n = 2v_n,\quad (v_n,v_n) + 2\alpha_n^2 
 = 2\alpha_n; \label{Aut:nn} \\
\alpha_i v_j = 0,  \quad (v_i,v_j) + 2\alpha_i\alpha_j = 0,
\quad i,j<n,\,i\neq j;  \label{Aut:ij} \\
\alpha_i v_n = 0,  \quad (v_n,v_i) + 2\alpha_n\alpha_i = 0,
\quad i<n;  \label{Aut:in} \\
\alpha_n v_j = v_j,\quad (v_n,v_j) + 2\alpha_n\alpha_j 
 = \alpha_j,\quad j<n. \label{Aut:nj}
\end{gather}
Since $v_j$ are not all zero for $j=1,\ldots,n-1$ 
(otherwise, $\varphi$ is denegerate), we conclude by~\eqref{Aut:nj} 
that $\alpha_n = 1$. By~\eqref{Aut:nn}, we get $v_n = 0$ 
(hence, $\varphi(e_n) = e_n$) and so, by~\eqref{Aut:in}, $\alpha_i = 0$ for all $i<n$.
From the whole system, only the equalities
$(v_i,v_i) = 1$ and $(v_i,v_j) = 0$ for $i,j<n$, $i\neq j$, remain. 
They mean exactly that $\varphi|_{J_n}$ is an orthogonal matrix.

Now, consider $d\in\Der(I_n)$ and set $d(e_i) = w_i + \gamma_i e_n$ 
for $w_i\in J_n$ and $\gamma_i\in F$. Moreover, 
$w_i = w_{1i}e_1 + \ldots + w_{n-1\,i}e_{n-1}$. 
By the Leibniz rule,
$$
d(e_i e_j)
 = d(e_i) e_j + e_i d(e_j)
 = (w_i+\gamma_i e_n)e_j + e_i(w_j+\gamma_j e_n).
$$
Again, we derive the system of equations,
\begin{gather}
\gamma_i e_i = w_n,\quad 2w_{ii} = \gamma_n,\quad i<n;\label{Der:ii} \\
w_n = 2w_n,\quad 4\gamma_n = 2\gamma_n; \label{Der:nn} \\
\gamma_i e_j = 0, \quad w_{ij} + w_{ji} = 0,\quad i,j<n,\,i\neq j;\label{Der:ij}\\
w_{in} + 2\gamma_i = 0,\quad i<n;  \label{Der:in} \\
\gamma_n e_j = 0,\quad w_{jn} + \gamma_j = 0,\quad j<n. \label{Der:nj}
\end{gather}
The equalities~\eqref{Der:nn} imply that $d(e_n) = 0$.
By~\eqref{Der:in}, $\gamma_i = 0$ for every $i<n$.
Therefore, only the relations 
$2w_{ii} = w_{ij}+w_{ji} = 0$ for $i,j<n$, $i\neq j$, remain. 
It means that $d|_{J_n}$ is a~skew-symmetric matrix.
\hfill $\square$

\section{Rota---Baxter operators}

Given an algebra $A$ and a scalar $\lambda\in F$, where $F$ is a~ground field,
a~linear operator $R\colon A\rightarrow A$ is called a Rota---Baxter operator
(RB-operator, for short) of weight~$\lambda$ if the identity
\begin{equation}\label{RB}
R(x)R(y) = R( R(x)y + xR(y) + \lambda xy )
\end{equation}
holds for all $x,y\in A$.

{\bf Proposition} \cite{GuoMonograph}.
Let an algebra $A$ split as a vector space
into a direct sum of two subalgebras $A_1$ and $A_2$.
An operator $P$ defined as
\begin{equation}\label{Split}
P(a_1 + a_2) = -\lambda a_2,\quad a_1\in A_1,\ a_2\in A_2,
\end{equation}
is RB-operator of weight~$\lambda$ on~$A$.

\newpage

Let us call an RB-operator from Proposition as
splitting RB-operator with subalgebras $A_1,A_2$.
It is known that $R$ is splitting if and only if $R^2+\lambda R = 0$~\cite{GuoMonograph}.

Let $A$ be an algebra and $R$ be a Rota---Baxter operator of weight~$\lambda$ on~$A$. 
It is well-known that $\phi(R)=-R-\lambda\id$ is again a Rota---Baxter operator of weight~$\lambda$ on~$A$.

Let $R$ be a Rota---Baxter operator of weight~$\lambda$ on $I_n$.
We denote $R(e_i) = v_i + \alpha_i e_n$, where $v_i\in J_n$ and 
$\alpha_i\in F$.
On the one hand, we have
\begin{equation} \label{RB-LHS}
R(e_i)R(e_j)
 = (v_i+\alpha_i e_n)(v_j+\alpha_j e_n)
 = \alpha_i v_j + ( (v_i,v_j) + 2\alpha_i \alpha_j )e_n.
\end{equation}
On the other hand,
\begin{multline} \label{RB-RHS}
R( R(e_i)e_j + e_iR(e_j) + \lambda e_i e_j) \\
 = R( \alpha_i(1+\delta_{n,j})e_j + v_{ji}I_{j<n}e_n  
 + \delta_{n,i}v_j + (2\alpha_j\delta_{n,i}+v_{ij}I_{i<n})e_n
 + \lambda e_i e_j ).
\end{multline}

Since~\eqref{RB-LHS} and~\eqref{RB-RHS} have to be equal,
we write down the following system of relations:
\begin{gather} 
(v_{ij}+v_{ji}+\lambda\delta_{i,j})v_n = 0,  \quad i,j<n,\label{RB-ij1} \\
(v_i,v_j) + \alpha_i \alpha_j 
 = (v_{ij}+v_{ji}+\lambda\delta_{i,j})\alpha_n, \quad i,j<n;\label{RB-ij2} \\
(\alpha_i+v_{in})v_n = 0, \quad i<n, \label{RB-in1} \\ 
(v_i,v_n) = v_{in}\alpha_n, \quad i<n, \label{RB-in2} \\
\lambda v_j + \sum\limits_{k=1}^{n-1}v_{kj}v_k + (v_{jn}+2\alpha_j)v_n = 0,
 \quad j<n, \label{RB-nj1} \\ 
\sum\limits_{k=1}^{n-1}\alpha_k v_{kj} + \lambda \alpha_j 
 + (v_{jn}+\alpha_j)\alpha_n = (v_n,v_j), \quad j<n, \label{RB-nj2} \\ 
\allowdisplaybreaks
(3\alpha_n+2\lambda)v_n + \sum\limits_{k=1}^{n-1}v_{kn}v_k = 0,\label{RB-nn1}\\
\sum\limits_{k=1}^{n-1}\alpha_k v_{kn} + 2\alpha_n(\alpha_n+\lambda)
 = (v_n,v_n). \label{RB-nn2}
\end{gather}

{\sc Case 1}: $v_n\neq0$.
Then we have some kind of skew-symmetricity by~\eqref{RB-ij1} and~\eqref{RB-in1}
\begin{equation} \label{Skew-sym}
v_{ij}+v_{ji}+\lambda\delta_{i,j} = 0, \quad
v_{in} = -\alpha_i, \quad i,j<n.
\end{equation}
From the first equalities, we find $v_{ii} = -\lambda/2$, $i<n$. 
The second ones imply that $\alpha_i\neq0$ for some $1\leq i\leq n-1$, since $v_n\neq0$.
The equalities~\eqref{RB-ij2} and~\eqref{RB-in2} may be rewritten as follows,
\begin{equation} \label{Scalar-prod}
(v_i,v_j) = -\alpha_i \alpha_j, \quad
(v_i,v_n) = -\alpha_i \alpha_n, \quad i,j<n.
\end{equation}

Now, we project~\eqref{RB-nn1} on $e_j$ for $j<n$ with the help of~\eqref{Skew-sym}:
$$
0 = (3\alpha_n + 2\lambda)v_{jn} + \sum\limits_{k=1}^{n-1}v_{kn}v_{jk}
  = -(3\alpha_n + 2\lambda)\alpha_j - \sum\limits_{k=1}^{n-1}\alpha_k v_{jk}.
$$
By the first equality of~\eqref{Skew-sym}, we rewrite it as follows,
\begin{equation} \label{RB-nn1-new}
0 = -(3\alpha_n + 2\lambda)\alpha_j 
  + (\lambda/2)\alpha_j
  + \sum\limits_{k=1,\,k\neq j}^{n-1}\alpha_k v_{kj}.
\end{equation}
Applying~\eqref{RB-nj2}, we express
$$
(3\alpha_n + 3/2\lambda)\alpha_j
 = \sum\limits_{k=1,\,k\neq j}^{n-1}\alpha_k v_{kj}
 = -\alpha_j(\alpha_n+\lambda/2).
$$
Hence,
$\alpha_j(2\alpha_n+\lambda) = 0$ for all $j<n$.
Since not all $\alpha_j$ are zero, $\alpha_n = -\lambda/2$.

By the second relation of~\eqref{Skew-sym}, 
we convert~\eqref{RB-nn2} to the following one:
\begin{equation} \label{alpha-n-2}
\alpha_n^2 + \sum\limits_{k=1}^{n-1}\alpha_k^2 
 = \alpha_n(2\alpha_n+\lambda) 
 = 0.
\end{equation}

One can show that~\eqref{RB-nj1} and~\eqref{RB-nn1} via its form~\eqref{RB-nn1-new} are equivalent 
to the first and to the second relations of~\eqref{Scalar-prod} respectively. 

Now, introduce a matrix $A = (a_{ij})_{i,j=1}^n\in M_n(F)$ such that
$a_{ij} = \begin{cases}
v_{ij}, & i<n, \\
\alpha_j, & i=n.
\end{cases}$
Denote by $E$ the identity matrix from $M_n(F)$.
The conditions~\eqref{Skew-sym},~\eqref{Scalar-prod}, and~\eqref{alpha-n-2}
may be reinterpreted as follows,
$A = S -\lambda E/2$ for a skew-symmetric matrix~$S$ such that
$S^2 = (\lambda^2/4)E$ (via $A^TA=0$).
When $\lambda \neq 0$ and $n$ is odd, 
the last equality contradicts to the property that 
rank of a~skew-symmetric matrix is even.
Moreover,
$$
A^2 + \lambda A
 = (S + \lambda E/2)(S -\lambda E/2)
 = S^2 - (\lambda^2/4)E
 = 0.
$$

{\sc Case 2}: $v_n = 0$. By~\eqref{RB-nn2}, 
$\alpha_n = 0$ or $\alpha_n = -\lambda$. 
Up to $\phi$, we may assume that $\alpha_n = 0$.
Hence,~\eqref{RB-ij2} means that $A^T A = 0$. 
The remaining nontrivial equalities~\eqref{RB-nj1} and~\eqref{RB-nj2} 
joint give that $A^2 + \lambda A = 0$.

{\bf Example 4}~\cite{preLie}.
The following operator~$R$ on $I_n$ over~$\mathbb{C}$ for $n\geq3$ 
$$
R(e_1) = e_n + \dfrac{1}{\sqrt{2-n}}\sum\limits_{i=2}^{n-1}e_i,\quad
R(e_i)=0,\ 2\leq i\leq n
$$
is an RB-operator $R$ of weight~0 coming from Case~2.

{\bf Example 5}.
The following operator~$R$ on $I_4$ over~$\mathbb{C}$
$$
R(e_1) = -ie_3 - e_4,\quad
R(e_2) =   e_3 -ie_4,\quad
R(e_3) =  ie_1 - e_2,\quad
R(e_4) =   e_1+ ie_2
$$
is an RB-operator $R$ of weight~0 coming from Case~1.
The corresponding matrix~$A$ is skew-symmetric (but not skew-Hermitian!).

Now, we are able to describe RB-operators of nonzero weight~$\lambda$ on $I_n$.
Call a subspace~$W$ of $I_n$ a Lagrangian one, if 
$(a,b)_{\ext} = 0$ for all $a,b\in W$, where
$$
(v+\alpha e_n, u+\beta e_n)_{\ext} = (v,u) + \alpha\beta, 
 \quad v,u\in J_n,\ \alpha,\beta\in F,
$$
the inner product in the whole $I_n$.

Up to $\phi$, we assume that $e_n\in \ker(R)$ and so, 
$\ker(R+\lambda\id)$ is a Lagrangian subalgebra of $I_n$, since for all
$x = v + \alpha e_n$ and 
$y = u + \beta e_{n}$, where $v,u\in J_n$, $\alpha,\beta\in F$, 
we have 
$xy = \alpha y + (x,y)_{\ext}e_n$ by~\eqref{In:MultTable}.
If $x,y\in \ker(R+\lambda\id)$, then $(x,y)_{\ext} = 0$, otherwise 
$e_n\in \ker(R)\cap \ker(R+\lambda\id)$, a~contradiction.

Hence, every decomposition of $I_n$ into a direct vector space sum 
of two subalgebras has the form $I_n = W \oplus (U\oplus \Span\{e_n\})$, where
$W$ is a Lagrangian subalgebra of~$I_n$, $U$ is a~complement of~$W$ into~$J_n$.

{\bf Example 6}.
We have two decompositions
$I_2 = \Span\{e_1\pm ie_2\}\oplus\Span\{e_2\}$ over~$\mathbb{C}$. 
Actually, all nontrivial RB-operators of nonzero weight on $I_2$ over~$\mathbb{C}$
are exhausted by splitting ones corresponding to the stated above decompositions. 

In~\cite{AnBai}, it is written that there are 
no nontrivial RB-operators of weight~$-1$ on $I_2$ over~$\mathbb{C}$.
Thus, the authors of~\cite{AnBai} missed these two decompositions.

Now, we describe all RB-operators of weight~0 on $I_n$ from Case~2.
For this, we interpret both conditions $R^2 = 0$ and $A^TA = 0$.
The last one means that $\Imm(R)$ is a Lagrangian subalgebra of~$I_n$.
Therefore, we may again decompose $I_n = W \oplus (U\oplus \Span\{e_n\})$, 
where $W$ is a Lagrangian subalgebra of~$I_n$, $U$ is a~complement of~$W$ into~$J_n$. 
Here we define $R(e_n) = 0$, $R(W) = 0$, and $R$ somehow maps~$U$ to~$W$.

Therefore, we have proved the following result.

{\bf Theorem 2}.
Let $F$ be a field of characteristic not~2.
Then for every Rota---Baxter operator~$R$ of weight~$\lambda$ on~$I_n$, 
we have $R^2+\lambda R = 0$ and $A^TA=0$, where $A$ is a matrix of~$R$ 
in the basis~$e_1,\ldots,e_n$.

{\bf Corollary 1}.
Let $F$ be a field of characteristic not~2.
Then every Rota---Baxter operator~$R$ of nonzero weight~$\lambda$ on~$I_n$ 
corresponds to a direct decomposition of~$I_n$ into a~sum of two subalgebras, i.\,e. $R$ is splitting.

A field~$F$ is called totally real, if $\nu_1^2+\ldots+\nu_t^2=0$
implies $\nu_1=\ldots=\nu_t=0$.
A field~$F$ is called quadratically closed, 
if every quadratic equation over~$F$ has a solution in it.

{\bf Corollary 2}.
Let $F$ be a field of characteristic not~2.

a) If $F$ is totally real field, then there are only trivial RB-operators on $I_n$.

b) If $F$ is quadratically closed, then there are 
nontrivial RB-operators of nonzero weight on $I_n$ for $n\geq2$ and 
nontrivial RB-operators of weight~0 on $I_n$ for $n\geq3$.

{\sc Proof}.
Given an RB-operator~$R$ of weight~$\lambda$ on $I_n$, 
either $R$ has weight~0 and its matrix~$A$ is skew-symmetric one 
such that $A^2 = -A^TA = 0$ (Case 1), or $R(e_n) = 0$, $R^2=-\lambda R$, 
and $A^TA=0$ (Case 2). In both cases, if $R$ is nontrivial, 
one has to solve quadratic equations 
$v_{i1}^2 + \ldots + v_{i\,n-1}^2 = 0$, $i=1,\ldots,n-1$ and 
additionally~\eqref{alpha-n-2} for $\lambda = \alpha_n = 0$ (from Case~2). 
Hence, we have the restrictions from b) and c).
On the other hand, we may construct nontrivial RB-operators 
of weight~$\lambda$ on $I_n$ by Example~4 and trivial extension of Example~6.
\hfill $\square$

The next definition was given in~\cite{Spectrum}.
Given an algebra $A$, denote the set of all RB-operators 
of weight~$\lambda$ on~$A$ by $\RB_\lambda(A)$.
The Rota---Baxter $\lambda$-index of $A$ is defined as follows,
$$
\rb_\lambda(A) = \min\{n\in\mathbb{N}\mid \mbox{for all }
 R\in \RB_\lambda(A)\ \mbox{exists } k:\ R^k(R+\lambda\id)^{n-k} = 0\}.
$$
If such number is undefined, put $\rb_\lambda(A) = \infty$.

{\bf Corollary 3}.
Let $F$ be a~field of characteristic not~2.
Then $\rb_\lambda(I_n)\leq2$.
If $F$ is totally real, then $\rb_\lambda(I_n)=1$.
If $F$ is quadratically closed and $n\geq3$, then $\rb_\lambda(I_n)=2$.

\section{Comments}

It is easy to show that one may join a unit element to every pre-Lie algebra~$A$ over a field~$F$ in an usual way, i.\,e. we define $A^\# = A\oplus F1$
and put $(a+\alpha1)(b+\beta1) = ab + \alpha b+\beta a + \alpha\beta1$ for $a,b\in A$ and $\alpha,\beta\in F$. It is again a pre-Lie algebra.

{\bf Remark 1}.
Note that all described in the work automorphisms, derivations, and Rota---Baxter operators may be extended onto $I_n^\#$ as follows.
If $\varphi\in\Aut(I_n)$, we put $\varphi(1) = 1$.
For $d\in\Der(I_n)$ and an RB-operator $R$ on $A$, we define $d(1) = R(1) = 0$.

{\bf Remark 2}.
If we consider $I_n$ under the anticommutator $a\circ b = (ab+ba)/2$,
we get a~commutative, not power-associative, flexible algebra $I_n^{(+)}$, 
which is simple if $\mathrm{char}\,F = 0$ or $\mathrm{char}\,F > 3$.
When $\mathrm{char}\,F = 3$ and $F$ is quadratically closed, then $I_n$ for $n\geq2$ is not simple.

{\bf Remark 3}.
Based on~$I_n$, one may construct a simple infinite-dimensional pre-Lie algebra.
Let us consider a vector space $I_\infty$ with a linear basis $e_0,e_1,\ldots$
under the product
$$
e_0 \cdot e_0 = 2e_0,\quad
e_0 \cdot e_j =  e_j,\quad
e_j \cdot e_j =  e_0,\quad j\geq 1.
$$
Then $I_\infty$ is a simple pre-Lie algebra.

\section*{Acknowledgments}

The author is grateful to Mikhail Pirozhkov, 
with whom we studied automorphisms and derivations of~$I_n$. 
Th author is also grateful to the anonymous reviewer.
The research is supported by Russian Science Foundation (project 21-11-00286).

\section*{Conflict of interests}

The author states that there is no conflict of interest.

\medskip
\noindent Vsevolod Gubarev \\
Sobolev Institute of Mathematics \\
Acad. Koptyug ave. 4, 630090 Novosibirsk, Russia \\
Novosibirsk State University \\
Pirogova str. 2, 630090 Novosibirsk, Russia \\
e-mail: wsewolod89@gmail.com \\
ORCID: 0000-0002-7839-5714
\end{document}